\documentclass[12pt]{article}

\makeatletter
\newfont{\footsc}{cmcsc10 at 8truept}
\newfont{\footbf}{cmbx10 at 8truept}
\newfont{\footrm}{cmr10 at 10truept}
\renewcommand{\ps@plain}{}
\makeatother 

\title{An old problem of Erd\H os: a graph without two cycles of the same length }

\author{ Chunhui Lai \thanks{ Project supported by the NSF of Fujian (2015J01018; 2020J01795; 2021J02048), Fujian Provincial Training Foundation for "Bai-Quan-Wan Talents Engineering", Fujian Key Laboratory of Granular Computing and Applications, Institute of Meteorological Big Data-Digital Fujian and  Fujian Key Laboratory of Data Science and Statistics  (Minnan Normal University) , China
}\\
\small School of Mathematics and Statistics, Minnan Normal University,
\\
\small Zhangzhou, Fujian 363000, P. R. of CHINA\\
\small \texttt{ laich2011@msn.cn; laichunhui@mnnu.edu.cn }\\
\small MR Subject Classifications: 05C38, 05C35\\
\small Key words: graph, cycle, number of edges}
\date{}
\begin{document}
\maketitle

\begin{abstract}
In 1975, P. Erd\H os  proposed the problem of determining the
maximum number $f(n)$ of edges in a  graph on $n$ vertices in
which any two cycles are of different
 lengths. Let $f^{\ast}(n)$ be the maximum number of edges in a simple graph on $n$ vertices
in which any two cycles are of different
 lengths. Let $M_n$  be the set of simple graphs on $n$ vertices and $f^{\ast}(n)$ edges in which any two cycles are of different lengths. Let  $mc(n)$ be the maximum cycle length for all  $G \in M_n$.
 In this paper, it is proved that  for $n$ sufficiently large, $mc(n)\leq \frac{15}{16}n.$
  \par

  We make the following conjecture:
  \par
  \bigskip
  \noindent{\bf Conjecture.} $$\lim_{n \rightarrow \infty} {mc(n)\over  n}= 0.$$
 \end{abstract}

\section{Introduction}
Let $f(n)$ be the maximum number of edges in a graph on $n$
vertices in which no two cycles have the same length. In 1975,
P. Erd\H os  raised the problem of determining $f(n)$ (see [1],
p.247, Problem 11). Let $f^{\ast}(n)$ be the maximum number of edges in a simple graph on $n$ vertices
in which any two cycles are of different
 lengths. Let $M_n$  be the set of simple graphs on $n$ vertices and $f^{\ast}(n)$ edges in which any two cycles are of different lengths. Let  $mc(n)$ be the maximum cycle length for all  $G \in M_n$. Let  $sc(n)$ be the second-largest cycle length for all  $G \in M_n$. Let  $tc(n)$ be the third-largest cycle length for all  $G \in M_n$. A natural question is what is the numbers of $mc(n)$, $sc(n)$, $tc(n).$ Let  $mcn(n)$ be the maximum cycle numbers for all  $G \in M_n$. A natural question is what is the numbers of $mcn(n).$ Let  $b(n)$ be the maximum $2$-connected block numbers for all  $G \in M_n$. A natural question is what is the numbers of $b(n).$
  Shi[23] proved that
   \par
\bigskip
  \noindent{\bf Theorem 1 (Shi [23]).}
$$f(n)\geq n+
 [(\sqrt {8n-23} +1)/2]$$ for $n\geq 3$ and $f(n)= f^{\ast}(n-1)+ 3$ for $n\geq 3$.
 \par
 Lai[8] proved that
 \par
\bigskip
  \noindent{\bf Theorem 2 (Lai [8]).}
 For $n\geq
e^{2m}(2m+3)/4$, $$f(n)<
 n-2+\sqrt {n ln (4n/(2m+3)) +2n} +log_2 (n+6).$$
\par
 Chen, Lehel, Jacobson and Shreve[3] gave a quick proof of this result.
 \par
Jia[6], Lai[7,8,9,10,11,12,13,14,15],
 Shi[23,24,25,26,27,28], Shi, Tang, Tang, Gong, Xu[29], Shi, Xu, Chen, Wang[30] obtained some additional related results.
 \par
 Lai[16] proved that
 \par
\bigskip
  \noindent{\bf Theorem 3 (Lai [16]).}
  $$\liminf_{n \rightarrow \infty} {f(n)-n\over \sqrt n}\geq \sqrt {2 +
\frac{40}{99}}.$$
and  Lai[9] conjectured that
\par
\bigskip
  \noindent{\bf conjecture 4 (Lai [9]).}
$$\liminf_{n \rightarrow \infty} {f(n)-n\over \sqrt n} \leq \sqrt {3}.$$
 Boros, Caro, F\"{u}redi and Yuster[2]
 proved that
  \par
\bigskip
  \noindent{\bf Theorem 5 (Boros, Caro, F\"{u}redi and Yuster[2]).}
  $$f(n)\leq n+1.98\sqrt{n}(1+o(1)).$$
\par

  Let $f_2(n)$ be the maximum number of edges in a $2$-connected
 graph on $n$ vertices in which no two cycles have the same length.
 \par
 In 1988, Shi[23] proved that\par
\par
\bigskip
  \noindent{\bf Theorem 6 (Shi[23]).}
 For every integer $n\geq3$, $f_{2}(n)\leq
n+[\frac{1}{2}(\sqrt{8n-15}-3)]$.
\par
 In 1998, G. Chen, J. Lehel, M. S. Jacobson, and W. E. Shreve [3]
 proved that
 \par
\par
\bigskip
  \noindent{\bf Theorem 7 (Chen, Lehel, Jacobson and Shreve [3]).}
 $f_{2}(n)\geq n+\sqrt{n/2}-o(\sqrt{n})$\par

  In 2001, E. Boros, Y. Caro, Z. F\"uredi and R. Yuster [2] improved this lower bound
significantly:\par
\par
\bigskip
  \noindent{\bf Theorem 8 (Boros, Caro, F\"{u}redi and Yuster[2]).}
 $f_{2}(n)\geq n+
\sqrt{n}-O(n^{\frac{9}{20}})$.
\par
and  conjectured that
\par
\par
\bigskip
  \noindent{\bf Conjecture 9 (Boros, Caro, F\"{u}redi and Yuster[2]).}
 \ $\lim \frac{f_{2} (n)-n}{\sqrt{n}}=1$.
\par
It is easy to see that this Conjecture  implies the (difficult)
upper bound in the Erd\H os Tur\'an Theorem [4][5](see [2]).
\par

\par
Markstr\"om [22] raised the problem
\par
\bigskip
  \noindent{\bf Problem 10 (Markstr\"om [22]).}
  Determining
 the maximum number of edges in a Hamiltonian
graph on $n$ vertices with no repeated cycle lengths.
 \par
 \par
Let $g(n)$ be the maximum number edges in an $n$-vertex, Hamiltonian graph with no repeated
cycle lengths.  J. Lee, C. Timmons [18] proved the following.
\par
\par
\bigskip
  \noindent{\bf Theorem 11 (J. Lee, C. Timmons [18]).}
 If $q$ is a power of a prime and $n = q^2 + q + 1,$ then
$$g(n)\geq n + \sqrt{n - 3/4} - 3/2$$

\par
A simple counting argument shows that  $g(n) < n + \sqrt{2n} +1$.

\par
Let $MH_n$  be the set of Hamiltonian graphs on $n$ vertices and $g(n)$  edges in which any two cycles are of different lengths. Let  $mcn_{H}(n)$ be the maximum cycle numbers for all  $G \in MH_n$. A natural question is what is the numbers of $mcn_{H}(n).$
\par

\par
J. Ma, T. Yang [21]  proved that
\par
\par
\bigskip
  \noindent{\bf Theorem 12 (Ma, Yang [21]).}
Any $n$-vertex $2$-connected graph with no two cycles of the same length contains at
most $n+\sqrt{n}+o(\sqrt{n})$ edges.
\par
Let $f_{2}(n, k)$ be the maximum number of edges in a graph $G$ on $n$ verticesin which no two cycles have the same length and $G$ which consists of $k$ $2$-connected blocks.
A natural question is what is the  maximum number of edges $f_{2}(n, k)$. It is clearly that $f_{2}(n, 1)=f_{2}(n)$.
\par
By theorem 5, it is clearly that
$$f_{2}(n, k) \leq f(n)\leq n+1.98\sqrt{n}(1+o(1)).$$
 \par
H. Lin, M. Zhai,Y. Zhao [19]  proved that
\par
\par
\bigskip
  \noindent{\bf Theorem 13 (Lin, Zhai,Zhao [19]).}
Let $G$ be a graph of order $n \geq 26.$ If $\rho(G) \geq \rho(K^{+}_{1,n-1}),$ then $G$ contains
two cycles of the same length unless $G \cong K^{+}_{1,n-1}.$
\par
and asked the following problem.
\par
\par
\bigskip
  \noindent{\bf Problem 14 (Lin, Zhai,Zhao [19]).}
What is the maximum spectral radius among all $2$-connected $n$-vertex
graphs without two cycles of the same length?
\par
Y. Shi [27]proved that
\par
\bigskip
  \noindent{\bf Theorem 15 (Shi [27]).}
$$b(n)\leq  [(\sqrt {8n+1} -5)/2]+1$$
\par
C. Lai [7]  proved that
\par
\par
\bigskip
  \noindent{\bf Theorem 16 (Lai [7]).}
$mc(n)\leq n-1$  for $n\geq \sum_{i=1}^{71}i - 8\times 18.$
\par
\par
 Survey papers on this problem can be found in Tian[31],
 Zhang[32], Lai and Liu[17].
\par
The progress of all 50 problems in [1] can be found in Locke[20].
Let $v(G)$ denote the number of vertices, and $\varepsilon(G)$ denote the number of edges. In this paper, it is proved that
   \par
\bigskip
  \noindent{\bf Theorem 17.}
 For $n$ sufficiently large, $$ mc(n)\leq \frac{15}{16}n.$$
  \par

 \section{Proof of the theorem 17}
  {\bf Proof.} If $ mc(n)> \frac{15}{16}n,$ for $n$ sufficiently large,
  then there is a simple graph $G$ on $n$ vertices and $f^{\ast}(n)$ edges in which any two cycles are of different
 lengths, the maximum cycle length of $G$ is $mc(n).$ Let $G_1$ be the block contain the cycle with length $mc(n).$
   It is clear that $v(G_1) > \frac{15}{16}n.$
 By the result of Ma and Yang [21],
 $\varepsilon (G_1) \leq v(G_1) + \sqrt{v(G_1)} + o(\sqrt{v(G_1)} ).$
 By the result of Boros, Caro, F\"uredi and Yuster [2],
 $\varepsilon (G) \leq v(G_1) + \sqrt{v(G_1)} + o(\sqrt{v(G_1)} ) + V(G)- V(G_1) + 1 + 1.98\sqrt{V(G)- V(G_1) + 1}(1+o(1)) \leq n+1+\sqrt{n}+o(\sqrt{n})+1.98\sqrt{\frac{1}{16}n}(1+o(1))\leq n+\frac{3}{2}\sqrt{n},$ for $n$ sufficiently large.
 By the result of Shi [23] and Lai [16], $\varepsilon (G) = f^{\ast}(n) = f(n+1)-3 > n + (\sqrt {2 +
\frac{40}{99}}-o(1))\sqrt{n},$ for $n$ sufficiently large. Note that  $\varepsilon (G) \leq n+\frac{3}{2}\sqrt{n},$
  this contradiction completes the proof.
 \vskip 0.2in
 \par
 It is clear that $mcn(n)\leq mc(n) - 2.$
  \par
 By theorem 3, it is clearly that
 $$mcn(n)\geq   \sqrt {2 + \frac{40}{99}}\sqrt{n}(1-o(1)).$$
  \par
  We make the following conjecture:
  \par
  \bigskip
  \noindent{\bf Conjecture.} $$\lim_{n \rightarrow \infty} {mc(n)\over  n}= 0.$$
 \par

\par
 \section*{Acknowledgment}
 The author thanks Professor Yair Caro and Raphael Yuster for sending me references[2].
 The author thanks Professor Yaojun Chen for sending me references[24].
 The author thanks Professor Tianchi Yang for sending me references[21].
 The author would like to thank Professor Yair Caro for his advice.
 \par

\end{document}